\documentclass[twocolumn]{IEEEtran}
\usepackage{amsmath,amssymb,amsthm,epsfig,color,empheq,graphicx,graphics,balance}
\usepackage{enumerate,url,wasysym,epstopdf,enumitem,array}
\usepackage{xcolor}
\usepackage{amsfonts}
\usepackage{adjustbox}
\usepackage{multirow}
\usepackage{cite}
\usepackage{accents}

\usepackage{algorithm2e} 
\usepackage{float} 

\DeclareMathOperator{\diag}{dg}

\newtheorem{proposition}{Proposition}

\newcommand \bzero{\mathbf{0}}
\newcommand \bone{\mathbf{1}}

\newcommand \bg{\mathbf{g}}

\newcommand \bp{\mathbf{p}}
\newcommand \bq{\mathbf{q}}

\newcommand \bv{\mathbf{v}}

\newcommand \by{\mathbf{y}}
\newcommand \bz{\mathbf{z}}

\newcommand \bI{\mathbf{I}}

\newcommand \bR{\mathbf{R}}

\newcommand \bU{\mathbf{U}}

\newcommand \bX{\mathbf{X}}

\newcommand \balpha{\boldsymbol{\alpha}}

\newcommand \bdelta{\boldsymbol{\delta}}

\newcommand \bsigma{\boldsymbol{\sigma}}

\newcommand \bLambda{\mathbf{\Lambda}}

\newcommand \mcB{\mathcal{B}}

\newcommand \mcO{\mathcal{O}}

\newcommand \mcZ{\mathcal{Z}}

\newcommand \tbv{\tilde{\mathbf{v}}}

\newcommand \tbz{\tilde{\mathbf{z}}}

\newcommand \tbalpha{\tilde{\boldsymbol{\alpha}}}

\newcommand \tbdelta{\tilde{\boldsymbol{\delta}}}

\newcommand \hbq{\hat{\mathbf{q}}}

\newcommand \bbq{\bar{\mathbf{q}}}

\newcommand \bbv{\bar{\mathbf{v}}}

\begin{document}

\title{Scalable Optimal Design of Incremental Volt/VAR Control using Deep Neural Networks}


\author{
    Sarthak Gupta,~\IEEEmembership{Graduate Student Member,~IEEE,}
    Ali~Mehrizi-Sani,~\IEEEmembership{Senior Member,~IEEE,}
    Spyros~Chatzivasileiadis,~\IEEEmembership{Senior Member,~IEEE,}
    and
    Vassilis Kekatos,~\IEEEmembership{Senior Member,~IEEE}

}	
	

	


\maketitle

\begin{abstract}
Volt/VAR control rules facilitate the autonomous operation of distributed energy resources (DER) to regulate voltage in power distribution grids. According to non-incremental control rules, such as the one mandated by the IEEE Standard 1547, the reactive power setpoint of each DER is computed as a piecewise-linear curve of the local voltage. However, the slopes of such curves are upper-bounded to ensure stability. On the other hand, incremental rules add a memory term into the setpoint update, rendering them universally stable. They can thus attain enhanced steady-state voltage profiles. Optimal rule design (ORD) for incremental rules can be formulated as a bilevel program. We put forth a scalable solution by reformulating ORD as training a deep neural network (DNN). This DNN emulates the Volt/VAR dynamics for incremental rules derived as iterations of proximal gradient descent (PGD). The rule parameters appear as DNN weights. To reduce the DNN depth, we leverage Nesterov’s accelerated PGD iterations. Analytical findings and numerical tests corroborate that the proposed ORD solution can be neatly adapted to single/multi-phase feeders.
\end{abstract}
	


\section{Introduction}
Local Volt/VAR (Volt-Ampere Reactive) control facilitates voltage regulation on distribution grids by providing reactive power compensation from DERs equipped with smart inverters. Different from centralized control schemes which incur large computational and communication burden, local rules decide DER setpoints based on local measurements. Volt/VAR control  rules can be categorized into \emph{non-incremental} and \emph{incremental} ones. The former compute DER reactive power setpoints based on local voltage readings. The IEEE Standard 1547.8 prescribes such non-incremental control rules as piecewise-linear functions of voltage~\cite{IEEE1547}. On the other hand, incremental Volt/VAR rules compute the \emph{change} in VAR setpoints as a function of voltage~\cite{LQD14,FZC15,7361761,VKZG16,Guido16}.

The existing literature on designing Volt/VAR control rules can be classified into \emph{stability}- and \emph{optimality}-centric works. Stability-centric works study the effect of Volt/VAR rules as a closed-loop dynamical system, which may be rendered unstable under steep slopes of non-incremental rules~\cite{FCL13,9091863}. In fact, to ensure stability, non-incremental rules may have to compromise on the quality of their steady-state voltage profile~\cite{VKZG16,9091863}. Incremental rules however do not experience stability limitations and can thus achieve improved voltage profiles compared to their non-incremental counterparts. Nonetheless, such improvements may come at the expense of longer settling times of the associated Volt/VAR dynamics~\cite{9091863}. 

Optimality-centric works focus on designing stable control rules to minimize a voltage regulation objective. To this end, optimization-based strategies have been employed to design affine non-incremental rules using heuristics~\cite{6601722,6727491,8365842}. Two of our recent works in~\cite{MGCK23} and~\cite{GCK22} have addressed the problem of optimally designing the slope, deadband, saturation, and reference voltage. Reference~\cite{MGCK23} performs ORD via a bilevel optimization applicable to single-phase feeders. Reference~\cite{GCK22} proposes DNN-based digital twins that emulate Volt/VAR dynamics, and reformulates ORD as a DNN training task for single-/multi-phase feeders. 

This letter deals with optimally selecting the shape of incremental Volt/VAR control rules, with contributions on three fronts: \emph{c1)} Although this \emph{optimal rule design} (ORD) task can be posed as a mixed-integer nonlinear optimization program, it does not scale well with the numbers of DERs, nodes, and grid loading scenarios. To address this challenge, the genuine idea here is to reformulate ORD as a deep-learning task and judiciously adapt the fast software modules widely available for training deep neural networks (DNNs). We have put forth a similar approach for designing non-incremental control rules in~\cite{GCK22}. However, migrating from non-incremental to incremental rules is non-trivial due to the different curve shapes, stability, and settling time properties. \emph{c2)} To further expedite ORD for incremental rules, we suggest implementing accelerated Nesterov-type variants of the rules to yield a shallower DNN emulator. \emph{c3)} We also establish the convergence of incremental rules on multiphase feeders. 

Recently, reference~\cite{Baosen22} deals with the optimal design of incremental rules. It uses DNNs with a single hidden layer to model piecewise-linear functions and formulates ORD as a reinforcement learning task. While~\cite{Baosen22} also utilizes DNNs to design incremental rules, we delineate from it in several ways. Reference~\cite{Baosen22} focuses on voltage control during transient dynamics, whereas this work aims at ORD to drive steady-state voltages closer to unity and over different grid loading scenarios. Reference~\cite{Baosen22} utilizes a DNN to model the piecewise-linear mapping of the rule. In contrast, this work develops a DNN-based digital twin that emulates end-to-end Volt/VAR dynamics. Lastly, we provide stability and convergence analysis for single- and multiphase feeders alike, whereas \cite{Baosen22} applies only to single-phase feeders.

The rest of this letter is organized as follows. Section~\ref{sec:model} models the feeder and discusses non-incremental and incremental Volt/VAR control rules. Section~\ref{sec:1pDNN} formulates DNN-based digital twins for Volt/VAR dynamics of incremental rules, and their accelerated version. It also presents ORD for single-phase feeders as a deep learning task. Section~\ref{sec:3pDNN} extends the ORD process to multiphase feeders. The incremental rules are then benchmarked against non-incremental rules from~\cite{GCK22} using tests on real-world data, in Section~\ref{sec:tests}. The letter is concluded in Section~\ref{sec:conclusions}.

\section{Volt/VAR Control Rules}\label{sec:model}
\allowdisplaybreaks
Consider a radial feeder serving $N$ buses equipped with DERs, indexed by $n$. Let $(\bq^{\ell},\bq)$ collect reactive loads and generations at all nodes. Vectors $(\bp,\bv)$ collect the net active power injections and voltage magnitudes at all nodes. The impact of $\bq$ on $\bv$ can be approximately captured using the linearized grid model~\cite{GCK22}
\begin{equation}\label{eq:ldf}
\bv\simeq \bX\bq+\tbv
\end{equation}
where $\tbv:=\bR\bp-\bX\bq^{\ell}+v_0\bone$ models the underlying \emph{grid conditions}, and $v_0$ is the substation voltage. Vector $\tbv$ represents the impact of non-controlled quantities $(\bp,\bq^{\ell})$ on voltages. Matrices $(\bR,\bX)$ depend on the feeder topology. For single-phase feeders, they are symmetric positive definite with positive entries~\cite{TJKT-SG21}. For multiphase feeders, they are non-symmetric and have positive and negative entries~\cite{VKZG16,GCK22}. 

\begin{figure}[t]
	\centering
	\includegraphics[scale=0.62]{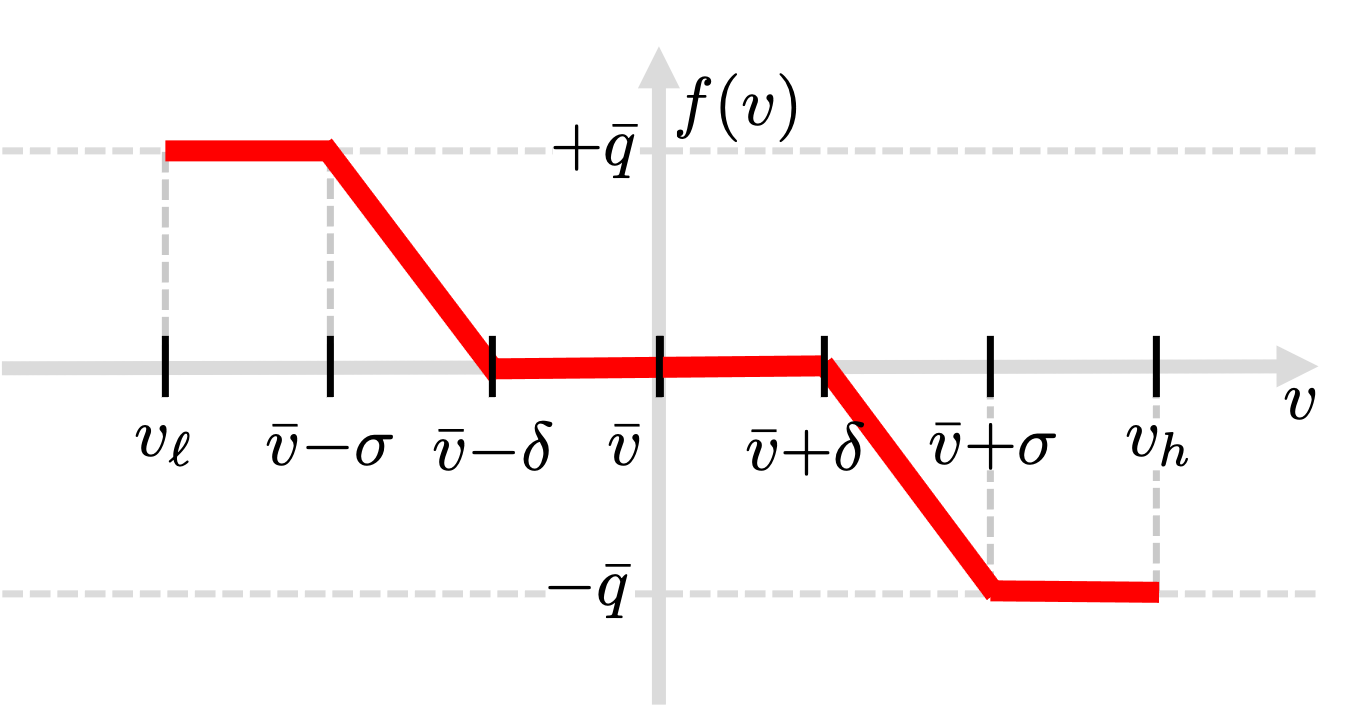}
	\caption{Non-incremental Volt/VAR control rule provisioned by the IEEE Std. 1547 for the interconnection of DERs~\cite{IEEE1547}.}
	\label{fig:curve}
\end{figure}

Vector $\bq$ in \eqref{eq:ldf} carries the reactive injections by DERs we would like to control. Per the non-incremental rules of the IEEE Std. 1547~\cite{IEEE1547}, DER setpoints are decided based on the Volt/VAR curve of Fig.~\ref{fig:curve}, which is parameterized by $(\bar{v},\delta,\sigma,\bar{q})$. The standard further constrains these parameters within a polytopic feasible set~\cite{IEEE1547,MGCK23}. The negative slope of the linear segment of the curve in Fig.~\ref{fig:curve} can be expressed as
\[\alpha:=\frac{\bar{q}}{\sigma-\delta}.\]
The interaction of Volt/VAR rules with the feeder gives rise to nonlinear dynamics. These dynamics are stable if $\|\diag(\balpha)\bX\|_2<1$, where $\diag(\balpha)$ is a diagonal matrix carrying the rule slopes over all buses on its diagonal~\cite{FCL13}. The equilibrium setpoints for DERs cannot be expressed in closed form. However, they coincide with the minimizer of the convex optimization problem~\cite{FCL13}
\begin{equation}\label{eq:inner}
\min_{-\bbq\leq \bq\leq\bbq}\frac{1}{2}\bq^\top \bX\bq+\bq^\top(\tbv-\bbv)+\frac{1}{2}\bq^\top\diag^{-1}(\balpha)\bq+\bdelta^\top|\bq|
\end{equation}
where $|\bq|$ applies the absolute value on $\bq$ entrywise. Problem~\eqref{eq:inner} depends on rule parameters $(\bar{v},\delta,\alpha,\bar{q})$ across all buses, collected in the $4N$-long vector $\bz:=(\bbv,\bdelta,\balpha,\bbq)$. We denote by $\bq_\bz(\tbv)$ the equilibrium setpoints, and by
\begin{equation}\label{eq:veq}
\bv_\bz(\tbv)=\bX\bq_\bz(\tbv)+\tbv
\end{equation} 
the related equilibrium voltages reached by Volt/VAR rules parameterized by $\bz$ under grid conditions $\tbv$.

\emph{Optimal rule design (ORD)} can be stated as the task of selecting $\bz$ to bring equilibrium voltages $\bv_\bz(\tbv)$ close to unity. To cater to diverse conditions,  the utility may sample loading scenarios $\{\tbv_s\}_{s=1}^S$ for the next hour, and find $\bz$ as
\begin{align}\label{eq:ord}\tag{ORD}
\bz^*\in \arg\min_{\bz}~&~F(\bz):=\frac{1}{S}\sum_{s=1}^S\|\bv_\bz(\tbv_s)-\bone\|_2^2\\
\textrm{subject to}~&~\eqref{eq:veq}~\text{and}~\bz\in\mcZ.\nonumber
\end{align}
Once found, the customized rules $\bz^*$ are sent to DERs to operate autonomously over the next hour. Note that $\bv_\bz(\tbv_s)$ depends on $\bz$ because the equilibrium setpoints $\bq_\bz(\bv_s)$ in \eqref{eq:veq} are the minimizers of problem \eqref{eq:inner}, which is parameterized by $\bz$. When solving \eqref{eq:ord} for non-incremental rules, the feasible set $\mcZ$ consists of the polytopic constraints imposed on $\bz$ by the IEEE Std. 1547 as well as additional constraints on $\balpha$ to ensure $\|\diag(\balpha)\bX\|_2<1$; see \cite{MGCK23}. Therefore, the feasible set $\mcZ$ can be quite confined. This can lead to less desirable voltage profiles; that is, higher objective values $F(\bz^*)$.

The aforesaid issue can be addressed by replacing the non-incremental Volt/VAR rules of IEEE Std. 1547 by incremental ones as suggested in~\cite{LQD14,FZC15,7361761,VKZG16,Guido16}. Incremental rules express the \emph{change} rather than the actual value in setpoints as a function of voltage. One option for incremental rules is to implement a proximal gradient descent (PGD) algorithm solving \eqref{eq:inner} as proposed in~\cite{VKZG16}. In this case, the control rule coincides with the PGD iterations, which are implemented by DERs in a decentralized fashion. Using incremental rules, set $\mcZ$ is enlarged as now we only need to ensure 
\begin{align*}
&\bz\geq \bzero\\
&0.95\cdot\bone\leq\bbv\leq 1.05\cdot\bone
\end{align*}
and that $\bbq$ are within the reactive power ratings of the DERs. 

The PGD algorithm is an extension of gradient descent to handle constraints and non-differentiable costs~\cite{VKZG16}. At iteration $t$, PGD proceeds with two steps: \emph{s1)} It first computes the gradient of the first two terms of $F(\bz)$, that is $\bX\bq^t+\tbv-\bbv=\bv^t-\bbv$. Here $\bq^t$ is the latest estimate of the minimizer of \eqref{eq:inner}; \emph{s2)} PGD then updates $\bq^{t+1}$ as the minimizer of
\begin{equation}\label{eq:pgds2}
\min_{-\bbq\leq \bq\leq \bbq}~\frac{1}{2}\bq^\top\diag^{-1}(\balpha)\bq +\bdelta^\top|\bq|+\frac{1}{2\mu}\|\bq-(\bv^t-\bbv)\|_2^2
\end{equation}
for a step size $\mu>0$. The last problem involves the last two terms in the cost of \eqref{eq:inner} regularized by the Euclidean distance of $\bq$ to the gradient $(\bv^t-\bbv)$ computed in step \emph{s1)}. 

Converting PGD to control rules, step \emph{s1)} is performed by the physics of the feeder when injecting $\bq^t$ and measuring the local voltage deviations. Step~\emph{s2)} is run by each DER independently as \eqref{eq:pgds2} is separable across buses. Using the subdifferential, solving \eqref{eq:pgds2} provides the update~\cite{VKZG16}
\begin{subequations}\label{eq:inc}
 \begin{align}
     y_n^t&=\tilde{\alpha}_n\cdot\left(q_n^t-\mu(v_n^t-\bar{v}_n)\right)\label{eq:inc:y}\\
     q_n^{t+1}&=g_n\left(y_n^t\right)\label{eq:inc:g}
 \end{align}
\end{subequations}
where $g_n(y_n)$ is the \emph{proximal operator} 
\begin{align}\label{eq:S}
    g_n(y_n):= 
    \begin{cases}
    +\bar{q}_n&,~ y_n>\overline{q}_n+\mu \tilde{\delta}_n\\ 
    y_n-\mu \tilde{\delta}_n &,~\mu \tilde{\delta}_n<y_n \leq \overline{q}_
n+\mu \tilde{\delta}_n \\ 
    0 &,~-\mu \tilde{\delta}_n \leq y_n \leq \mu \tilde{\delta}_n\\ 
    y_n+\mu\tilde{\delta}_n &,~-\overline{q}_n-\mu \tilde{\delta}_n \leq y_n<-\mu \tilde{\delta}_n \\ 
   -\overline{q}_n&,~ y_n<-\overline{q}_n-\mu \tilde{\delta}_n.
    \end{cases}
\end{align}
and the new parameters $(\tilde{\alpha}_n,\tilde{\delta}_n)$ are defined as
\[\tilde{\alpha}_n:=\frac{1}{1+\mu/\alpha_n}\quad \text{and}\quad \tilde{\delta}_n:=\frac{\delta_n}{1+\mu/\alpha_n}.\]
The proximal operator is plotted in the top panel of Figure~\ref{fig:relus_inc}. Note that in \eqref{eq:inc}, rule parameters are transformed from representation $\bz=(\bbv,\bdelta,\balpha,\bbq)$ to representation $\tbz:=(\bbv,\tbdelta,\tbalpha,\bbq)$. This is without loss of generality as the transformation is a bijection, and so one can work exclusively with $\tbz$. The feasible set $\tilde{\mcZ}$ for $\tbz$ is similar to $\mcZ$ with the addition that $\tbalpha\leq \bone$.
As with non-incremental rules, the rules in~\eqref{eq:S} are driven by local data, but now $q_n^{t+1}$ depends on $(v_n^t,q_n^t)$, and not $v_n^t$ alone. Both types of rules solve \eqref{eq:inner}. Hence, they both converge to the same equilibrium. The advantage of incremental rules is that they are stable for all $\balpha$ as long as $\mu<2/\lambda_{\text{max}}(\bX)$; see~\cite{VKZG16}. It is worth stressing that $\bz$ does not have the same physical interpretation as in non-incremental rules (slopes, deadband, or saturation), though $\bz$ parameterizes \eqref{eq:inner} for both rules. 

\emph{Accelerated incremental rules.} Although PGD rules enlarge $\mcZ$, their settling times can be long. They reach an $\varepsilon$-optimal cost of \eqref{eq:inner} within $-\frac{2\log{\varepsilon}}{\log{2}}\kappa\left(\bX\right)$ iterations. Here $\kappa(\bX):=\lambda_{\max}(\bX)/\lambda_{\min}(\bX)$ is the condition number of $\bX$. References~\cite{VKZG16,SGC15} put forth \emph{accelerated} incremental rules based on accelerated PGD (APGD). These rules need $-\frac{2\log{\varepsilon}}{\log{2}} \sqrt{\kappa\left(\bX\right)}$ iterations to attain an $\varepsilon$-optimal cost, and take the form
\begin{subequations}\label{eq:acc_inc}
 \begin{align}
     \tilde{y}_n^t&:=\left(1+\beta_t\right)y_n^t-\beta_t y_n^{t-1}\label{eq:acc_inc:ytilde}\\
     q_n^{t+1}&:=g_n\left(\tilde{y}_n^t\right)
 \end{align}
\end{subequations}
where $\beta_t:=\frac{t-1}{t+2}$, while $y_n^t$ and $g_n(y_n)$ are as defined in~\eqref{eq:inc:y} and \eqref{eq:S}. Updates \eqref{eq:acc_inc} remain local, but introduce additional memory as $q_n^{t+1}$ depends on $(v_n^t,q_n^t)$ and $(v_n^{t-1},q_n^{t-1})$. 

\section{Deep Learning for Optimal Rule Design (ORD) in Single-Phase Feeders}\label{sec:1pDNN}

\begin{figure}[t]
	\centering
	\includegraphics[scale=0.6]{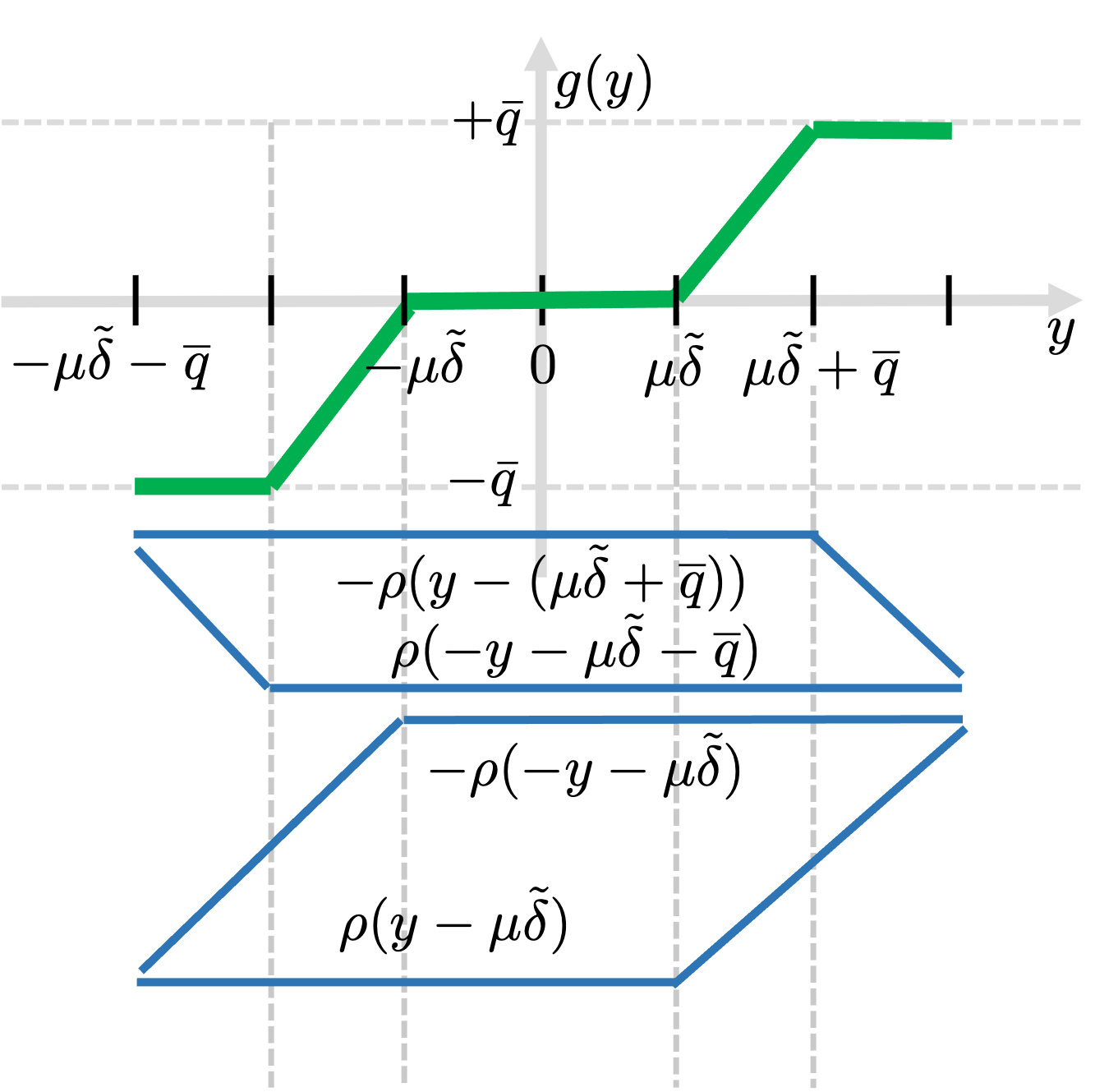}
	\caption{Proximal operator $g(y)$ expressed as a sum of four shifted rectified linear units (ReLUs).}
	\label{fig:relus_inc}
\end{figure}

Solving \eqref{eq:ord} is challenging as it is a nonconvex bilevel program. Although it can be modeled as a mixed-integer nonlinear program, such an approach does not scale well with the number of DERs and/or scenarios for non-incremental rules~\cite{GCK22}. Seeking a more scalable solution, we reformulate \eqref{eq:ord} as a deep learning task. The key idea is to design a DNN that emulates Volt/VAR dynamics under the control rule of \eqref{eq:inc}. To this end, note that $g_n(y_n)$ is a piecewise-linear function with four breakpoints~\cite{VKZG16}. Interestingly, this operator can be expressed as the superposition of four rectified linear units (ReLUs) as illustrated in Fig.~\ref{fig:relus_inc}, where ReLUs are denoted by $\rho(\cdot)$. The intercepts of the ReLUs depend linearly on $(\tilde{\delta}_n,\bar{q}_n)$. 

 \begin{figure}[t]
	\centering
	\includegraphics[scale=0.47]{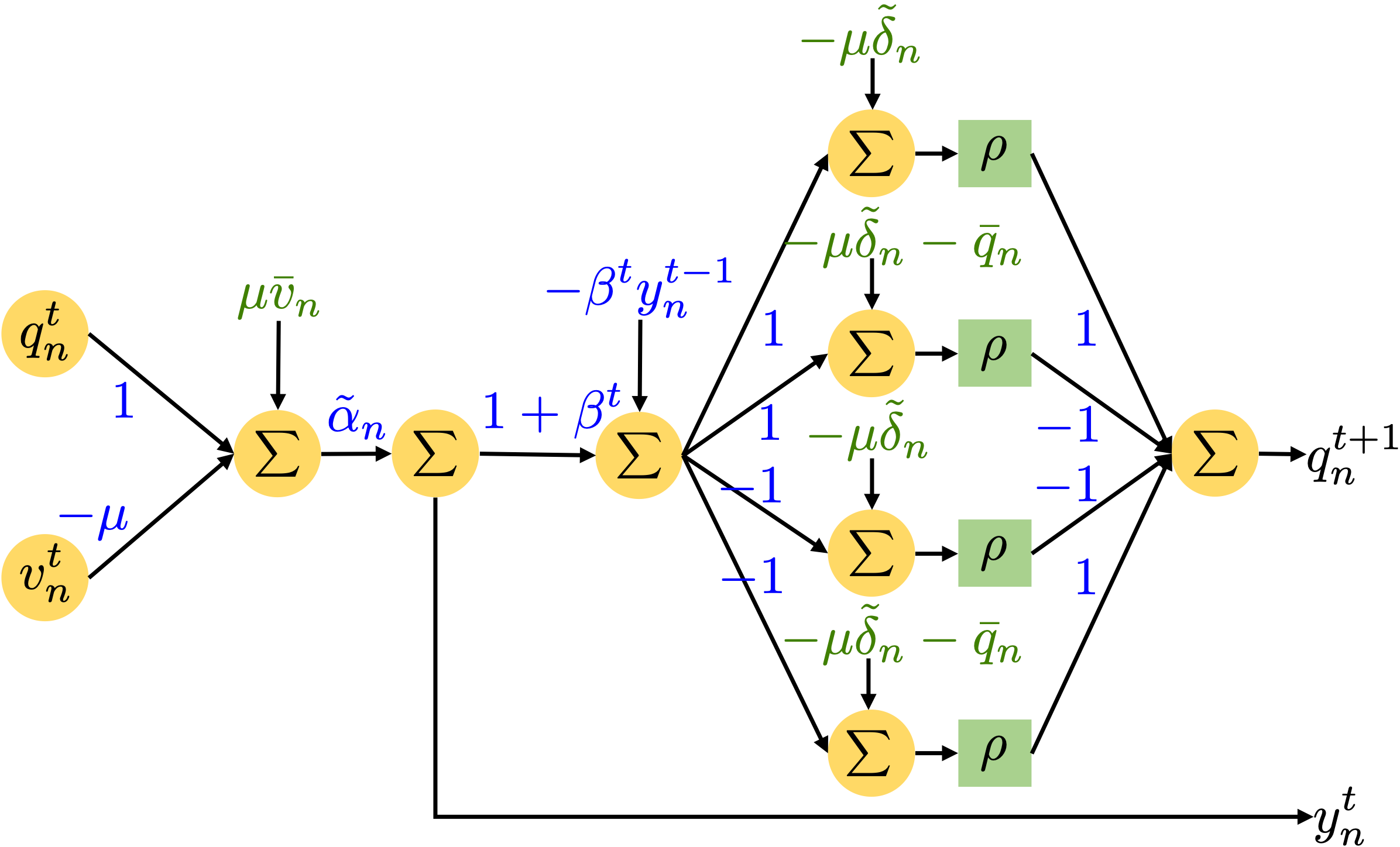}
	\caption{A DNN emulating the accelerated incremental rules of~\eqref{eq:acc_inc}. Plain incremental rules can be modeled by dropping the second layer (setting $\beta^t=0$) and ignoring output $y_n^t$.}
	\label{fig:nest_single}
\end{figure}

\begin{figure}[t]
	\centering
	\includegraphics[scale=0.46]{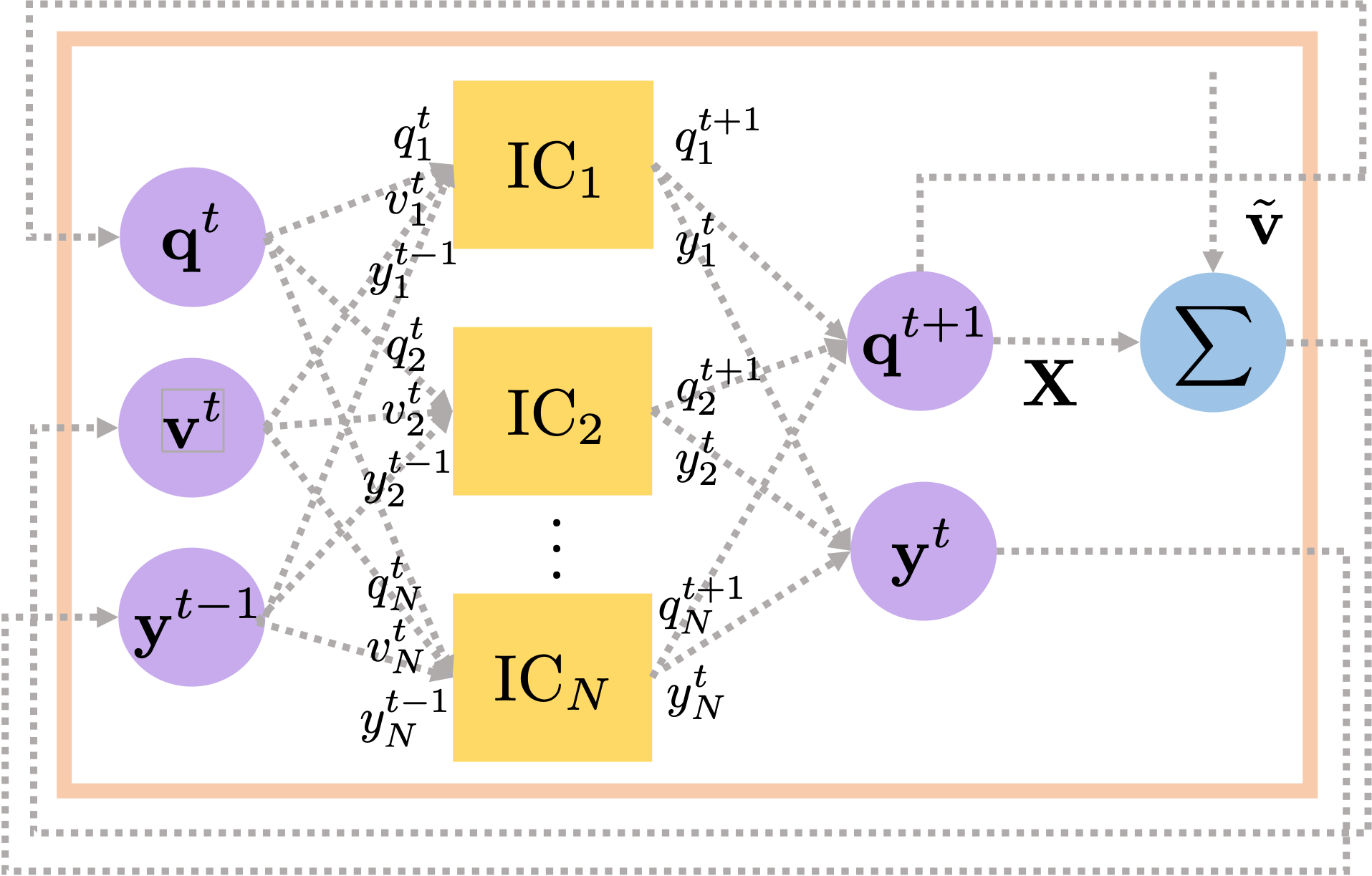}
	\caption{Recurrent neural network (RNN) implementation for accelerated incremental Volt/VAR control rules.}
	\label{fig:nest_RNN}
\end{figure}

Building on this, one APGD iteration for DER $n$ can be implemented by the 4-layer DNN in Fig.~\ref{fig:nest_single}, whose weights depend affinely on $(\bar{v}_n,\tilde{\delta}_n,\tilde{\alpha}_n,\bar{q}_n)$. This DNN takes $(q_n^t,v_n^t)$ as its input, and computes $(q_n^{t+1},y_n^t)$ at its output. It is termed $\text{IC}_n$ and will be used as a building block to emulate Volt/VAR dynamics. This is accomplished by the recursive neural network (RNN) shown in Fig.~\ref{fig:nest_RNN}. Here blocks $\text{IC}_n$ are arranged vertically to model the parallel operation of DERs. Their outputs $\bq^{t+1}$ are multiplied by $\bX$, and the new voltage is computed as $\bv^{t+1}=\bX\bq^{t+1}+\tbv$. This is repeated $T$ times. Thanks to the RNN structure, there is \emph{weight sharing}, so the number of DNN weights is $4N$ rather than $4NT$. 

The RNN takes a grid loading vector $\tbv_s$ as its input, the rule parameters $\tbz$ as weights, and computes the voltages $\bv^T_{\tbz}(\tbv_s)$ at time $T$ at its output. For the output $\bv^T_{\tbz}(\tbv_s)$ to approximate well equilibrium voltages, the depth $T$ can be chosen by the convergence rate of PGD as follows.

\begin{proposition}\label{pro:accdepth}
	For the DNN of Fig.~\ref{fig:nest_RNN} to ensure $\|\Phi\left(\tbv;\bz\right) -\bv^*(\bz)\|_2\leq \epsilon_1$ $\forall$ $\tbv$, its depth $T$ should satisfy
	\begin{equation}\label{eq:lemma:inc}
	T\geq \left(\frac{\kappa-1}{2}\right)\log\left(\frac{2\|\bX\|_2\|\hbq\|_2}{\epsilon_1}\right).
	\end{equation}
\end{proposition}

\emph{Proof:} From the control rule of~\eqref{eq:inc:g}, it follows that
\begin{align}
\|\bq^{t}-\bq^{*}\|_2&=\|\bg\left(\by^{t}\right)-\bg\left(\by^{*}\right)\|_2\nonumber\\
&\leq \|\by^{t}-\by^{*}\|_2\nonumber\\
&=\|\diag(\tilde{\balpha})(\bI-\mu\bX)\left(\bq^{t-1}-\bq^{*}\right)\|_2\nonumber\\
&\leq \|\diag(\tilde{\balpha})\|_2\cdot \|\bI-\mu\bX\|_2\cdot \|\bq^{t-1}-\bq^{*}\|_2\nonumber\\
&\leq\|\bI-\mu\bX\|_2\cdot \|\bq^{t-1}-\bq^{*}\|_2.\label{eq:pro:accdepth:bound}
\end{align}
The first inequality stems from the non-expansive property of the proximal operator $\bg$. The next equality follows from~\eqref{eq:inc:y}. The second inequality from the sub-multiplicative property of the spectral norm. The last inequality follows by the definition of spectral norm and because $\tilde{\alpha}_n\leq1$ for all $n$.

If $\|\bI-\mu\bX\|_2<1$, inequality~\eqref{eq:pro:accdepth:bound} implies that the dynamics in~\eqref{eq:inc} are a non-expansive mapping, and thus, are stable and converge to $\bq^{*}$. Condition $\|\bI-\mu\bX\|_2<1$ holds when $\mu<2/\lambda_{\max}(\bX)$. The norm $\|\bI-\mu\bX\|_2$ achieves its minimum of $\left(1-\frac{2}{\kappa+1}\right)$ when \[\mu_0:=\frac{2}{\lambda_{\max}(\bX)+\lambda_{\min}(\bX)}.\]
Plugging $\mu_0$ into~\eqref{eq:pro:accdepth:bound} and unfolding the dynamics over $t$ provides
\begin{align*}
\|\bq^t-\bq^*\|_2&\leq\left(1-\tfrac{2}{\kappa+1}\right)^t \|\bq^0-\bq^*\|_2\\
&\leq
2\left(1-\tfrac{2}{\kappa+1}\right)^t\|\hbq\|_2.
\end{align*}
For the voltage approximation error $	\|\bv^T-\bv^*\|_2 = \|\bX\left(\bq^T-\bq^*\right)\|_2$ at time $T$ to be smaller than $\epsilon_1$, we need
\begin{equation*}
\|\bv^T-\bv^*\|_2\leq 2\|\bX\|_2\cdot\|\hbq\|_2\cdot \left(1-\frac{2}{\kappa+1}\right)^T\leq \epsilon_1.
\end{equation*}
This can be achieved by selecting $T$ such that
\begin{align*}
T&\geq \frac{\log\left(\frac{2\|\bX\|_2\|\hbq\|_2}{\epsilon_1}\right)}{\log\left(1+\frac{2}{\kappa-1}\right)}\\
&\geq\left(\frac{\kappa-1}{2}\right)\log\frac{2\|\bX\|_2\|\hbq\|_2}{\epsilon_1}.
\end{align*}
where the last inequality follows from $\log(1+x)\leq x$. $\qed$

Plugging the values $\|\bX\|_2=0.463$ and $\kappa=848$ for the IEEE 37-bus feeder, $\|\hbq\|_2=0.1$, and $\epsilon_1=10^{-5}$ in \eqref{eq:lemma:inc}, yields $T\geq 2,892$ layers, which is relatively large. A key contributor to this large $T$ is the $\kappa$ term in \eqref{eq:lemma:inc}. This promulgates the adoption of accelerated rules~\eqref{eq:acc_inc}, which are known to have $\mcO(\sqrt{\kappa})$ dependence. Interestingly, during implementation, one does not need to fix $T$ to the above worst-case bounds. Leveraging dynamic computation graphs offered by Python libraries such as Pytorch, one may determine $T$ \emph{`on the fly'} depending on the convergence of $\bv^t$ between pairs of successive layers.

Since the RNN emulates Volt/VAR dynamics, it can surrogate $\bv_{\bz}(\tbv_s)$ in \eqref{eq:ord}. Then \eqref{eq:ord} can be posed as training a DNN over its weights $\tbz\in\tilde{\mcZ}$ or $\bz\in\mcZ$. Grid loading scenarios $\{\tbv_s\}_{s=1}^S$ are treated as features and equilibrium voltages $\bv_\bz(\tbv_s)$ as predictions that should be brought close to the target value of $\bone$ for scenarios $s$. The DNN can be trained using stochastic projected gradient descent (SPGD) as~\cite{GCK22} 
\begin{align}\label{eq:spgd}
\tbz^{i+1}&=\left[\tbz^{i}-\frac{\lambda}{2B}{\nabla_{\tbz^i}}\left(\sum_{s\in \mcB_i} \|\Phi(\tbv_s;\tbz)-\bone\|_2^2\right)\right]_{\tilde{\mcZ}}
\end{align}
where $\lambda>0$ is the learning rate; set $\mcB_i$ is a batch of $B$ scenarios; and $[\cdot]_{\tilde{\mcZ}}$ is the projection onto $\tilde{\mcZ}$. Since $\tilde{\mcZ}$ consists of simple box constraints, projection essentially means clipping the values to the box. Lastly $\nabla_{\tbz^i}(\cdot)$ represents the gradient with respect to $\tbz$ evaluated at $\tbz=\tbz^i$, and is calculated efficiently thanks to \emph{gradient back-propagation}. Although our DNN-based ORD assumed PGD-based rule, it may be applicable to other incremental rules too. 

A discussion about control rules and their DNN-based emulators is due. Recall that all three types of Volt/VAR control rules (non-incremental, incremental, and accelerated incremental) reach the same equilibrium voltages, if stable. The emulators aim at computing these equilibrium voltages. A natural question is whether the DNN emulator could implement a rule of a type different from the rule actually implemented on the feeder. This may be desirable to leverage the advantages of two types. Some caution is needed here. If the feeder implements non-incremental rules, but incremental rules converge faster to equilibrium voltages, it makes sense for the emulator to implement incremental rules. Of course, in this case, stability constraints on the non-incremental rules have to be enforced during DNN training. The reverse is not recommended: If the emulator implements non-incremental rules, its parameters $\bz$ should be constrained to be stable and that would be a restriction of the actual ORD problem. Finally, given the convergence advantage of accelerated incremental rules, they are always preferable over plain incremental rules for the DNN implementation. This showcases the utility of accelerated control rules even if they are not actually implemented on the feeder.


\section{Deep Learning for Optimal Rule Design (ORD) in Multiphase Feeders}\label{sec:3pDNN}
In multiphase feeders, matrix $\bX$ is non-symmetric and has both positive and negative entries. Therefore, the rule analysis and design of Section~\ref{sec:1pDNN} has to be revisited. For example, equilibrium setpoints cannot be found as the minimizers of an optimization problem as with \eqref{eq:inner}. Moreover, increasing $\bq$ does not mean that all voltages increase.

In multiphase feeders, the non-incremental rules of IEEE Std. 1547 remain stable as long as $\|\diag(\balpha)\bX\|_2<1$. This is the same condition as in the single-phase setup. How about the stability and equilibrium of incremental rules in multiphase feeders? Recall that for single-phase rules, incremental rules were obtained as the PGD iterations solving~\eqref{eq:inner}. Lacking an equivalent inner optimization for multiphase feeders precludes a similar approach here. Despite the incremental rules of~\eqref{eq:inc} do not correspond to PGD iterates anymore, they can still be shown to be stable for multiphase feeders. 

\begin{proposition}\label{prop:inc:3p:step}
Let $\bU\bLambda\bU^{\top}$ be the eigen-decomposition of matrix $\bX\bX^{\top}$. The incremental rules of~\eqref{eq:inc} are stable for multiphase feeders if their step size is selected as $\mu <\lambda_{\min}\left(\bLambda^{-1/2}\bU^{\top}\left(\bX+\bX^{\top}\right)\bU\bLambda^{-1/2}\right)$.
\end{proposition}

The claim follows readily by adopting the proof of Proposition~\ref{pro:accdepth}: If $\mu$ is selected as above, then $\|\bI-\mu\bX\|_2<1$ follows from~\cite[Prop. 6]{VKZG16}. Similar to the single-phase case, incremental rules in multiphase feeders allow us to enlarge the feasible set $\mcZ$ of rule parameters $\bz$. It is worth stressing that different from the single-phase setting, incremental and non-incremental rules do not converge to the same equilibrium on multiphase feeders.

The ORD task for multiphase feeders can also be formulated as a deep-learning task, with some modifications. Firstly, matrices $\bR$ and $\bX$ need to be altered. Secondly, the DNNs for multiphase feeders have $12N$ trainable parameters,  since each layer consists of $3N$ building modules corresponding to bus/phase (node) combinations. Lastly, the step size has to be selected per Proposition~\ref{prop:inc:3p:step}. Adopting the proof of Proposition~\ref{pro:accdepth}, we next find the minimum DNN depth in multiphase feeders. 

\begin{proposition}\label{pro:inc:3p}
Let the DNN of Fig.~\ref{fig:nest_RNN} implement the incremental rules of~\eqref{eq:inc} on multiphase feeders with $\mu$ selected per Proposition~\ref{prop:inc:3p:step}. The DNN depth $T$ ensuring voltage approximation error $\|\Phi\left(\tbv;\bz\right) -\bv^*(\bz)\|_2\leq \epsilon_1$ is
\begin{align*}
 T\geq \frac{ \log\frac{\epsilon_1}{2\|\bX\|_2\|\hbq\|_2}}{\log\|\bI-\mu\bX\|_2}.   
\end{align*}
\end{proposition}

We next numerically evaluate the proposed DNN-based ORD approach in single- and multiphase feeders, and contrast the performance of incremental control rules with that of non-incremental rules. 


\section{Numerical Tests}\label{sec:tests}
We benchmark the performance of DNN-based incremental rules against non-incremental rules from~\cite{GCK22} on single- and multiphase feeders. Real-world data were sourced from the Smart* project on April 2, 2011~\cite{Smartsolar}, as explained in~\cite{GCK22}. The DNNs were implemented and trained using Pytorch.

\begin{figure}
	\centering
	\includegraphics[scale=0.6]{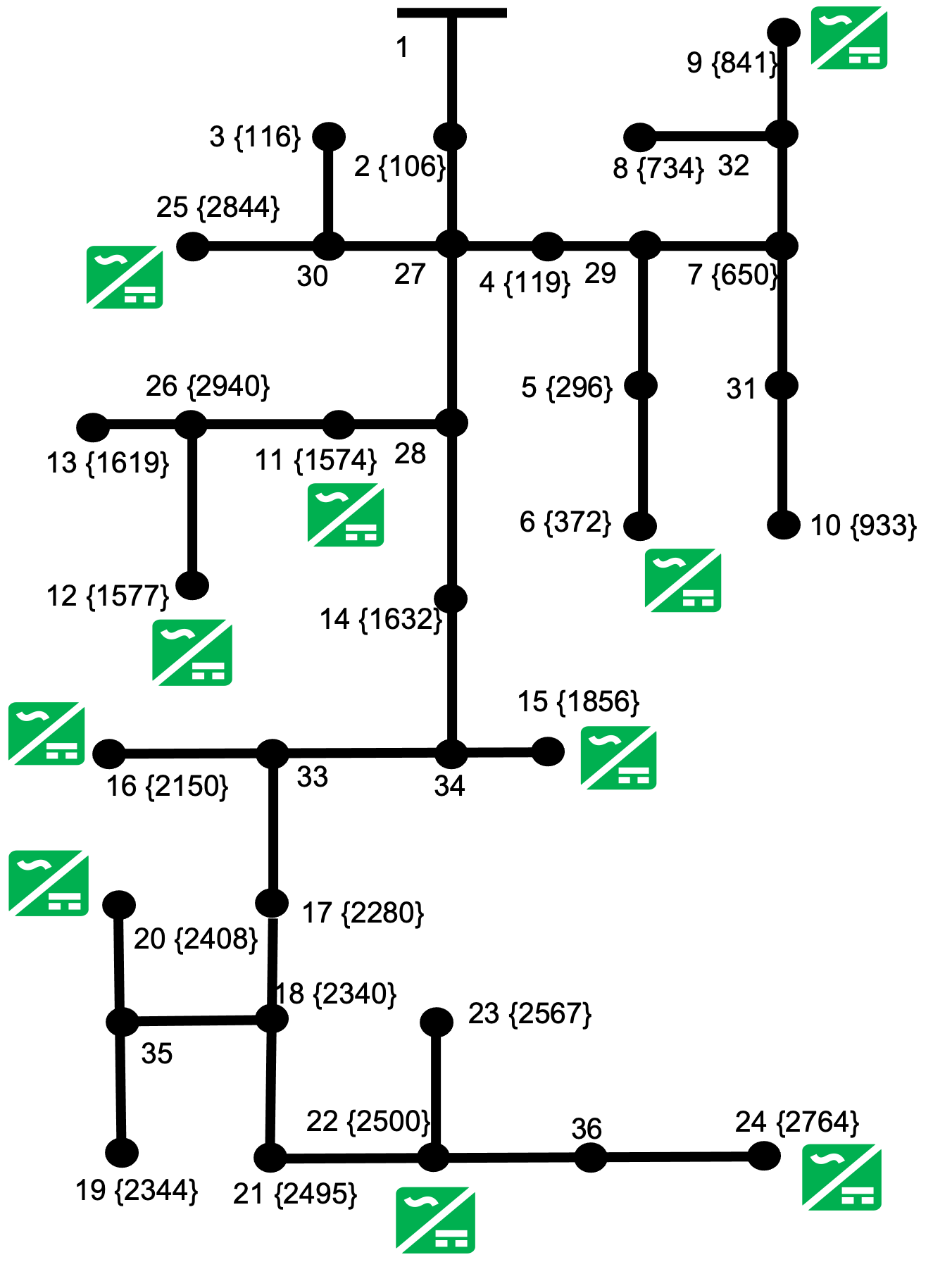}
	\caption{The IEEE 37-bus feeder converted to single-phase. Node numbering follows the format \texttt{node number \{panel ID\}}. DERs at buses $\{6,9,11,12,15,16,20,22,24,25\}$ provide reactive power control; the rest operate at unit power factor.}
	\label{fig:37bus}
\end{figure}

We first compare (non)-incremental rules, both designed via DNN training for the single-phase IEEE 37-bus feeder of Figure~\ref{fig:37bus}. Homes with IDs 20--369 were averaged 10 at a time and successively added as active
loads to buses 2--26 as shown in Fig. 6. Active generation from solar panels was also added, as per the mapping in
Fig. 6. Buses $\{6, 9, 11, 12, 15, 16, 20, 22, 24, 25\}$
were assumed to host DERs with Volt/VAR control customized per bus.

Incremental rules were simulated in their accelerated rendition. Both sets of rules were trained over $S=80$ scenarios and $200$ epochs with a learning rate of $0.001$, using the Adam optimizer, and setting $\mu=1$ for incremental rules. To ensure repeatability, the results were repeated across several time periods between 1--6~PM, and are compiled in Table~\ref{tab:nivsi}. Incremental rules obtained marginally lower objectives than non-incremental rules across all periods, with a somewhat significant difference for the 5~PM period. This behavior is explained because incremental rules allow for a larger set $\mcZ$. 

\begin{table}[t]
\centering\label{tab:nivsi}
\caption{Incremental vs. non-incremental Volt/VAR control rules on the single-phase IEEE 37-bus feeder}
\begin{adjustbox}{width=1\columnwidth}
\begin{tabular}{|c|c|cc|cc|}
\hline
\multirow{2}{*}{Time} & \multicolumn{1}{c|}{$\bq=\bzero$} & \multicolumn{2}{c|}{Non-incremental}      & \multicolumn{2}{c|}{Incremental}                       \\ \cline{2-6} 

                           & \multicolumn{1}{c|}{Obj. (p.u.)} & \multicolumn{1}{c|}{Time (s)} & Obj. (p.u.)                                 & \multicolumn{1}{c|}{Time (s)} & Obj. (p.u.) \\ \hline
 $1$ pm             & \multicolumn{1}{c|}{$3.01\cdot10^{-3}$} & \multicolumn{1}{c|}{$37.98$} & \multicolumn{1}{c|}{$3.68\cdot10^{-4}$}        & \multicolumn{1}{c|}{$39.39$} & $3.66\cdot10^{-4}$ \\ \hline
 $2$ pm             & \multicolumn{1}{c|}{$3.13\cdot10^{-3}$} & \multicolumn{1}{c|}{$42.93$} & \multicolumn{1}{c|}{$4.26\cdot10^{-4}$}        & \multicolumn{1}{c|}{$37.91$} & $4.25\cdot10^{-4}$ \\ \hline
    $3$ pm             & \multicolumn{1}{c|}{$4.24\cdot10^{-3}$} & \multicolumn{1}{c|}{$45.02$} & \multicolumn{1}{c|}{$8.59\cdot10^{-4}$}        & \multicolumn{1}{c|}{$34.97$} & $8.50\cdot10^{-4}$ \\ \hline
    $4$ pm             & \multicolumn{1}{c|}{$2.12\cdot10^{-3}$} & \multicolumn{1}{c|}{$48.30$} & \multicolumn{1}{c|}{$1.47\cdot10^{-4}$}        & \multicolumn{1}{c|}{$38.52$} & $1.48\cdot10^{-4}$ \\ \hline
    $5$ pm             & \multicolumn{1}{c|}{$8.53\cdot10^{-4}$} & \multicolumn{1}{c|}{$47 .37$} & \multicolumn{1}{c|}{$9.70\cdot10^{-5}$}        & \multicolumn{1}{c|}{$374.01$} & $6.90\cdot10^{-5}$ \\ \hline
\end{tabular}
\end{adjustbox}
\end{table}

\begin{figure}[t]
	\centering
	\includegraphics[scale=0.5]{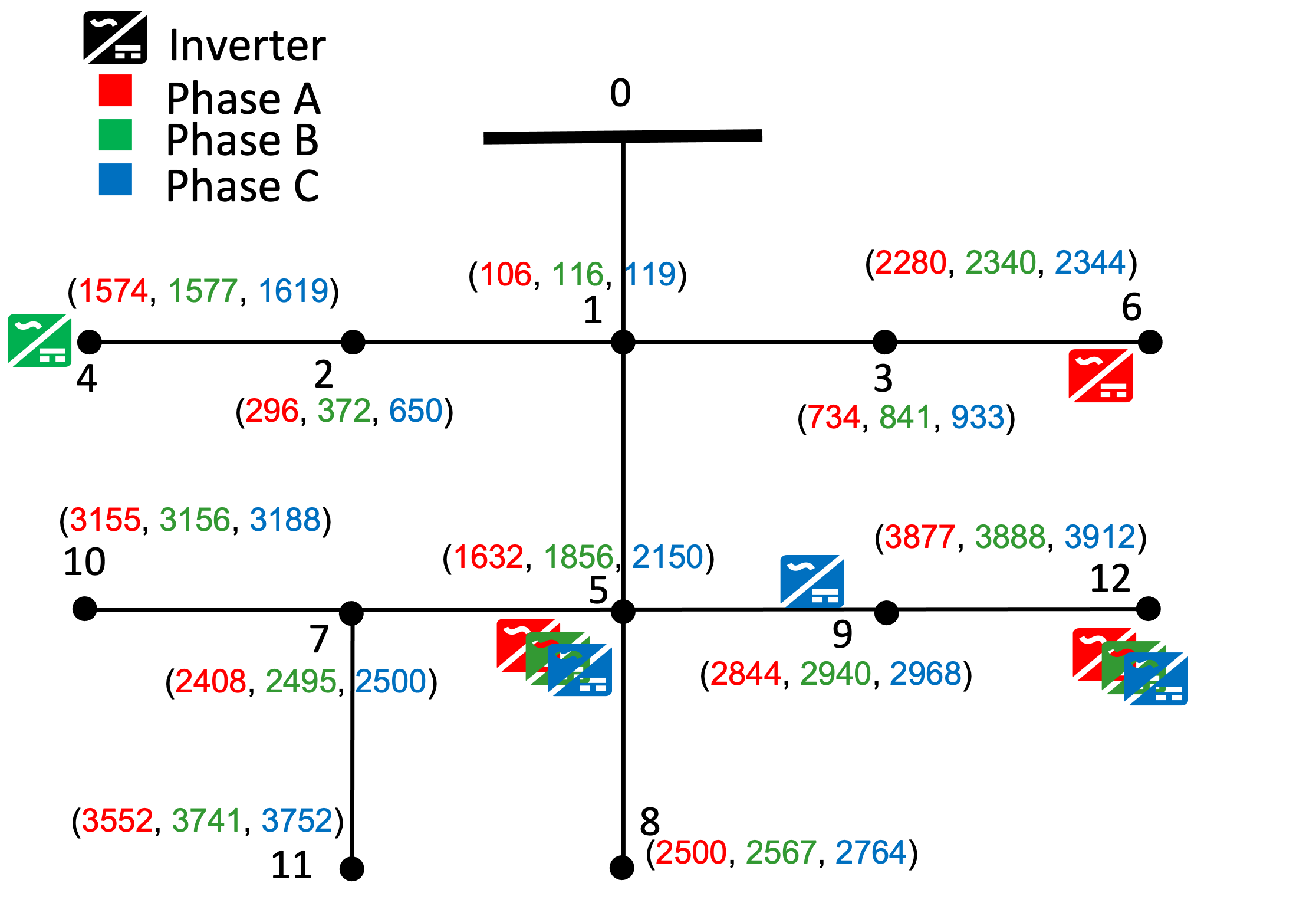}
	\caption{Multiphase IEEE 13-bus distribution feeder.}
	\label{fig:ieee13}
\end{figure}

DNN-based incremental control rules were also contrasted with their non-incremental ones on the multiphase IEEE 13-bus feeder, using the testing setup from~\cite{GCK22}.  Active loads were sampled $10$ at a time from homes with IDs $20$-$379$ and added to all three phases for the buses $1$-$12$. Figure~\ref{fig:ieee13} also shows the solar panel assignments shown in Fig~\ref{fig:ieee13} for solar generation. Lastly, nine DERs with inverters were added across phases and bus indices as shown in Fig~\ref{fig:ieee13}.

The learning rates for non-incremental and incremental DNNs were set as $0.1$ and $0.001$, respectively, with the design parameters $\bz:=(\bbv,\bdelta,\bsigma,\balpha)$ initialized to feasible values $(0.95,0.01,0.3,1.5)$. Table~\ref{tab:13_nivsi} compares the performance of the two rule categories over multiple periods for $S=80$. While incremental rules took longer times to train, they were successful in lowering the cost $F(\bz)$ by more than $50\%$, thus yielding improved voltage profiles across all periods.

\begin{table}[t]
\centering\label{tab:13_nivsi}
\caption{Incremental vs. non-incremental Volt/VAR control rules on the multiphase IEEE 13-bus feeder}
\begin{adjustbox}{width=1\columnwidth}
\begin{tabular}{|c|c|cc|cc|}
\hline
\multirow{2}{*}{Time} & \multicolumn{1}{c|}{$\bq=\bzero$} & \multicolumn{2}{c|}{Non-incremental}      & \multicolumn{2}{c|}{Incremental}                       \\ \cline{2-6} 

                           & \multicolumn{1}{c|}{Obj. (p.u.)} & \multicolumn{1}{c|}{Time (s)} & Obj. (p.u.)                                 & \multicolumn{1}{c|}{Time (s)} & Obj. (p.u.) \\ \hline
 $1$ pm             & \multicolumn{1}{c|}{$2.51\cdot10^{-3}$} & \multicolumn{1}{c|}{$64.65$} & \multicolumn{1}{c|}{$1.15\cdot10^{-3}$}        & \multicolumn{1}{c|}{$199.24$} & $4.11\cdot10^{-4}$ \\ \hline
 $2$ pm             & \multicolumn{1}{c|}{$1.48\cdot10^{-3}$} & \multicolumn{1}{c|}{$66.60$} & \multicolumn{1}{c|}{$6.89\cdot10^{-4}$}        & \multicolumn{1}{c|}{$209.92$} & $3.03\cdot10^{-4}$ \\ \hline
    $3$ pm             & \multicolumn{1}{c|}{$6.89\cdot10^{-4}$} & \multicolumn{1}{c|}{$74.68$} & \multicolumn{1}{c|}{$4.94\cdot10^{-4}$}        & \multicolumn{1}{c|}{$263.37$} & $2.16\cdot10^{-4}$ \\ \hline
    $4$ pm             & \multicolumn{1}{c|}{$8.03\cdot10^{-4}$} & \multicolumn{1}{c|}{$68.32$} & \multicolumn{1}{c|}{$5.26\cdot10^{-4}$}        & \multicolumn{1}{c|}{$126.81$} & $2.47\cdot10^{-4}$ \\ \hline
    $5$ pm             & \multicolumn{1}{c|}{$5.51\cdot10^{-4}$} & \multicolumn{1}{c|}{$62.58$} & \multicolumn{1}{c|}{$4.11\cdot10^{-4}$}        & \multicolumn{1}{c|}{$129.71$} & $1.95\cdot10^{-4}$ \\ \hline
\end{tabular}
\end{adjustbox}
\end{table}


\section{Conclusions}\label{sec:conclusions}
We have devised a DNN approach to optimally design incremental Volt/VAR control rules for single- and multi-phase feeders. The key idea is to construct a DNN that emulates end-to-end the associated Volt/VAR dynamics. The DNN takes grid conditions as the input, the rule parameters as weights, and outputs the associated equilibrium voltages. Leveraging the convergence rates of the related optimization algorithms, we have provided bounds on the minimum depth of the DNN emulator to approximate equilibrium voltages within the desired accuracy. We have also established the stability of incremental control rules for multiphase feeders. Numerical tests have demonstrated that the designed control rules attain improved voltage profiles compared to their non-incremental alternatives. The improvement was found to be starker for mutiphase feeders, wherein (non)-incremental rules do not reach the same equilibrium. Our findings motivate further research to possibly characterize the equilibria of control rules for multiphase feeders; the convergence of accelerated incremental rules for multiphase feeders; and to deal with chance-constrained formulations or ORD problems targeting phase imbalances.

\balance
\bibliographystyle{IEEEtran}
\bibliography{myabrv,kekatos,inverters}
\end{document}